\pdfoutput=1

\documentclass[12pt]{article}

\usepackage{amsfonts,amssymb,amsthm,amsmath,latexsym}
\usepackage{subeqnarray,eqsection,indent,eepicemu,url,cite,bm}
\usepackage{multirow}  
\usepackage{tensor}  

\usepackage[normalem]{ulem}

\date{December 7, 2016 \\[0mm] Corrigendum added: January 29, 2021 \\[3mm]
   Published in Lin. Alg. Appl. {\bf 520}, 242--259 (2017); \\[1.5mm]
   corrigendum in Lin. Alg. Appl. {\bf 613}, 393--396 (2021) }

\begin{document}

\title{\vspace*{-2.5cm}Total positivity of sums, Hadamard products \\
         and Hadamard powers: \\[3mm]
       Results and counterexamples \\[3mm]
       \hspace*{-8mm} {\bf WITH CORRIGENDUM AT END OF FILE}
      }

\author{
     {\small Shaun Fallat}  \\[-2mm]
     {\small\it Department of Mathematics and Statistics}       \\[-2mm]
     {\small\it University of Regina}       \\[-2mm]
     {\small\it Regina, Saskatchewan S4S 0A2}       \\[-2mm]
     {\small\it CANADA}       \\[-2mm]
     {\small\tt Shaun.Fallat@uregina.ca}   \\[-2mm]
     {\protect\makebox[5in]{\quad}}
     \\
     {\small Charles R.~Johnson}  \\[-2mm]
     {\small\it Department of Mathematics}       \\[-2mm]
     {\small\it College of William and Mary}     \\[-2mm]
     {\small\it Williamsburg, VA 23187}      \\[-2mm]
     {\small\it USA}      \\[-2mm]
     {\small\tt crjohn@wm.edu}   \\[-2mm]
     {\protect\makebox[5in]{\quad}}
     \\
     {\small Alan D.~Sokal\thanks{Also at Department of Mathematics,
        University College London, London WC1E 6BT, United Kingdom.}}  \\[-2mm]
     {\small\it Department of Physics}       \\[-2mm]
     {\small\it New York University}         \\[-2mm]
     {\small\it 4 Washington Place}          \\[-2mm]
     {\small\it New York, NY 10003}      \\[-2mm]
     {\small\it USA}      \\[-2mm]
     {\small\tt sokal@nyu.edu}               \\[-2mm]
     {\protect\makebox[5in]{\quad}}  
}

\maketitle
\thispagestyle{empty}   

\begin{abstract}
\noindent
We show that, for Hankel matrices,
total nonnegativity (resp.\ total positivity) of order $r$
is preserved by sum, Hadamard product, and Hadamard power
with real exponent $t \ge r-2$.
We give examples to show that our results are sharp
relative to matrix size and structure (general, symmetric or Hankel).
Some of these examples also resolve the Hadamard critical-exponent problem
for totally positive and totally nonnegative matrices.
\end{abstract}


\bigskip
\noindent
{\bf Mathematics Subject Classification (MSC 2010) codes:}
15B48 (Primary);
15A15, 15B05, 44A60 (Secondary).

\clearpage

\newtheorem{theorem}{Theorem}[section]
\newtheorem{proposition}[theorem]{Proposition}
\newtheorem{lemma}[theorem]{Lemma}
\newtheorem{corollary}[theorem]{Corollary}
\newtheorem{definition}[theorem]{Definition}
\newtheorem{conjecture}[theorem]{Conjecture}
\newtheorem{question}[theorem]{Question}
\newtheorem{problem}[theorem]{Problem}
\newtheorem{example}[theorem]{Example}

\renewcommand{\theenumi}{\alph{enumi}}
\renewcommand{\labelenumi}{(\theenumi)}
\def\eop{\hbox{\kern1pt\vrule height6pt width4pt
depth1pt\kern1pt}\medskip}
\def\prf{\par\noindent{\bf Proof.\enspace}\rm}
\def\rmk{\par\medskip\noindent{\bf Remark\enspace}\rm}

\newcommand{\be}{\begin{equation}}
\newcommand{\ee}{\end{equation}}
\newcommand{\<}{\langle}
\renewcommand{\>}{\rangle}
\newcommand{\widebar}{\overline}
\def\reff#1{(\protect\ref{#1})}
\def\spose#1{\hbox to 0pt{#1\hss}}
\def\ltapprox{\mathrel{\spose{\lower 3pt\hbox{$\mathchar"218$}}
    \raise 2.0pt\hbox{$\mathchar"13C$}}}
\def\gtapprox{\mathrel{\spose{\lower 3pt\hbox{$\mathchar"218$}}
    \raise 2.0pt\hbox{$\mathchar"13E$}}}
\def\textprime{${}^\prime$}
\def\proof{\par\medskip\noindent{\sc Proof.\ }}
\def\firstproof{\par\medskip\noindent{\sc First Proof.\ }}
\def\secondproof{\par\medskip\noindent{\sc Second Proof.\ }}
\renewcommand{\qed}{ $\square$ \bigskip}
\newcommand{\myendremark}{ $\blacksquare$ \bigskip}
\def\proofof#1{\bigskip\noindent{\sc Proof of #1.\ }}
\def\firstproofof#1{\bigskip\noindent{\sc First Proof of #1.\ }}
\def\secondproofof#1{\bigskip\noindent{\sc Second Proof of #1.\ }}
\def\half{ {1 \over 2} }
\def\third{ {1 \over 3} }
\def\twothird{ {2 \over 3} }
\def\smfrac#1#2{{\textstyle{#1\over #2}}}
\def\smhalf{ {\smfrac{1}{2}} }
\newcommand{\real}{\mathop{\rm Re}\nolimits}
\renewcommand{\Re}{\mathop{\rm Re}\nolimits}
\newcommand{\imag}{\mathop{\rm Im}\nolimits}
\renewcommand{\Im}{\mathop{\rm Im}\nolimits}
\newcommand{\sgn}{\mathop{\rm sgn}\nolimits}
\newcommand{\tr}{\mathop{\rm tr}\nolimits}
\newcommand{\supp}{\mathop{\rm supp}\nolimits}
\newcommand{\disc}{\mathop{\rm disc}\nolimits}
\newcommand{\diag}{\mathop{\rm diag}\nolimits}
\newcommand{\tridiag}{\mathop{\rm tridiag}\nolimits}
\newcommand{\perm}{\mathop{\rm perm}\nolimits}
\def\hboxscript#1{ {\hbox{\scriptsize\em #1}} }
\renewcommand{\emptyset}{\varnothing}

\newcommand{\restrict}{\upharpoonright}

\newcommand{\compinv}{{\langle -1 \rangle}}   

\newcommand{\scra}{{\mathcal{A}}}
\newcommand{\scrb}{{\mathcal{B}}}
\newcommand{\scrc}{{\mathcal{C}}}
\newcommand{\scrd}{{\mathcal{D}}}
\newcommand{\scre}{{\mathcal{E}}}
\newcommand{\scrf}{{\mathcal{F}}}
\newcommand{\scrg}{{\mathcal{G}}}
\newcommand{\scrh}{{\mathcal{H}}}
\newcommand{\scrk}{{\mathcal{K}}}
\newcommand{\scrl}{{\mathcal{L}}}
\newcommand{\scrm}{{\mathcal{M}}}
\newcommand{\scrn}{{\mathcal{N}}}
\newcommand{\scro}{{\mathcal{O}}}
\newcommand{\scrp}{{\mathcal{P}}}
\newcommand{\scrq}{{\mathcal{Q}}}
\newcommand{\scrr}{{\mathcal{R}}}
\newcommand{\scrs}{{\mathcal{S}}}
\newcommand{\scrt}{{\mathcal{T}}}
\newcommand{\scrv}{{\mathcal{V}}}
\newcommand{\scrw}{{\mathcal{W}}}
\newcommand{\scrz}{{\mathcal{Z}}}

\newcommand{\ahat}{{\widehat{a}}}
\newcommand{\Zhat}{{\widehat{Z}}}
\renewcommand{\k}{{\mathbf{k}}}
\newcommand{\n}{{\mathbf{n}}}
\newcommand{\vv}{{\mathbf{v}}}
\newcommand{\bv}{{\mathbf{v}}}
\newcommand{\w}{{\mathbf{w}}}
\newcommand{\x}{{\mathbf{x}}}
\newcommand{\cc}{{\mathbf{c}}}
\newcommand{\zero}{{\mathbf{0}}}
\newcommand{\one}{{\mathbf{1}}}
\newcommand{\bmm}{ {\bf m} }

\newcommand{\C}{{\mathbb C}}
\newcommand{\D}{{\mathbb D}}
\newcommand{\Z}{{\mathbb Z}}
\newcommand{\N}{{\mathbb N}}
\newcommand{\Q}{{\mathbb Q}}
\newcommand{\PP}{{\mathbb P}}
\newcommand{\R}{{\mathbb R}}
\newcommand{\RR}{{\mathbb R}}
\newcommand{\E}{{\mathbb E}}

\newcommand{\Sym}{{\mathfrak{S}}}
\newcommand{\SymB}{{\mathfrak{B}}}

\newcommand{\myle}{\preceq}
\newcommand{\myge}{\succeq}
\newcommand{\mygt}{\succ}

\newcommand{\B}{{\sf B}}
\newcommand{\OB}{{\sf OB}}
\newcommand{\OS}{{\sf OS}}
\newcommand{\OO}{{\sf O}}
\newcommand{\SP}{{\sf SP}}
\newcommand{\OSP}{{\sf OSP}}
\newcommand{\Eu}{{\sf Eu}}
\newcommand{\ERR}{{\sf ERR}}
\newcommand{\sfB}{{\sf B}}
\newcommand{\sfE}{{\sf E}}
\newcommand{\sfG}{{\sf G}}
\newcommand{\sfJ}{{\sf J}}
\newcommand{\sfP}{{\sf P}}
\newcommand{\sfQ}{{\sf Q}}
\newcommand{\sfS}{{\sf S}}
\newcommand{\sfT}{{\sf T}}
\newcommand{\sfW}{{\sf W}}
\newcommand{\sfMV}{{\sf MV}}
\newcommand{\AMV}{{\sf AMV}}
\newcommand{\BM}{{\sf BM}}

\newcommand{\emIB}{{\hbox{\em IB}}}
\newcommand{\emIP}{{\hbox{\em IP}}}
\newcommand{\emOB}{{\hbox{\em OB}}}
\newcommand{\emSC}{{\hbox{\em SC}}}

\newcommand{\stat}{{\rm stat}}
\newcommand{\cyc}{{\rm cyc}}
\newcommand{\Asc}{{\rm Asc}}
\newcommand{\asc}{{\rm asc}}
\newcommand{\Des}{{\rm Des}}
\newcommand{\des}{{\rm des}}
\newcommand{\Exc}{{\rm Exc}}
\newcommand{\exc}{{\rm exc}}
\newcommand{\Wex}{{\rm Wex}}
\newcommand{\wex}{{\rm wex}}
\newcommand{\Fix}{{\rm Fix}}
\newcommand{\fix}{{\rm fix}}
\newcommand{\lrmax}{{\rm lrmax}}
\newcommand{\rlmax}{{\rm rlmax}}
\newcommand{\Rec}{{\rm Rec}}
\newcommand{\rec}{{\rm rec}}
\newcommand{\Arec}{{\rm Arec}}
\newcommand{\arec}{{\rm arec}}
\newcommand{\ERec}{{\rm ERec}}
\newcommand{\erec}{{\rm erec}}
\newcommand{\EArec}{{\rm EArec}}
\newcommand{\earec}{{\rm earec}}
\newcommand{\Cpeak}{{\rm Cpeak}}
\newcommand{\cpeak}{{\rm cpeak}}
\newcommand{\Cval}{{\rm Cval}}
\newcommand{\cval}{{\rm cval}}
\newcommand{\Cdasc}{{\rm Cdasc}}
\newcommand{\cdasc}{{\rm cdasc}}
\newcommand{\Cddes}{{\rm Cddes}}
\newcommand{\cddes}{{\rm cddes}}
\newcommand{\Peak}{{\rm Peak}}
\newcommand{\peak}{{\rm peak}}
\newcommand{\Val}{{\rm Val}}
\newcommand{\val}{{\rm val}}
\newcommand{\Dasc}{{\rm Dasc}}
\newcommand{\dasc}{{\rm dasc}}
\newcommand{\Ddes}{{\rm Ddes}}
\newcommand{\ddes}{{\rm ddes}}
\newcommand{\inv}{{\rm inv}}
\newcommand{\maj}{{\rm maj}}
\newcommand{\rs}{{\rm rs}}
\newcommand{\cross}{{\rm cr}}
\newcommand{\crosshat}{{\widehat{\rm cr}}}
\newcommand{\nest}{{\rm ne}}
\newcommand{\lodd}{{\rm lodd}}
\newcommand{\leven}{{\rm leven}}
\newcommand{\sg}{{\rm sg}}
\newcommand{\bl}{{\rm bl}}
\newcommand{\tran}{{\rm tr}}
\newcommand{\area}{{\rm area}}
\newcommand{\ret}{{\rm ret}}
\newcommand{\peaks}{{\rm peaks}}
\newcommand{\hl}{{\rm hl}}
\newcommand{\sll}{{\rm sll}}
\newcommand{\negg}{{\rm neg}}

\newcommand{\ba}{{\bm{a}}}
\newcommand{\bahat}{{\widehat{\bm{a}}}}
\newcommand{\sfa}{{{\sf a}}}
\newcommand{\bb}{{\bm{b}}}
\newcommand{\bc}{{\bm{c}}}
\newcommand{\bff}{{\bm{f}}}
\newcommand{\bg}{{\bm{g}}}
\newcommand{\bu}{{\bm{u}}}
\newcommand{\bA}{{\bm{A}}}
\newcommand{\bB}{{\bm{B}}}
\newcommand{\bC}{{\bm{C}}}
\newcommand{\bE}{{\bm{E}}}
\newcommand{\bF}{{\bm{F}}}
\newcommand{\bI}{{\bm{I}}}
\newcommand{\bJ}{{\bm{J}}}
\newcommand{\bM}{{\bm{M}}}
\newcommand{\bN}{{\bm{N}}}
\newcommand{\bP}{{\bm{P}}}
\newcommand{\bQ}{{\bm{Q}}}
\newcommand{\bS}{{\bm{S}}}
\newcommand{\bT}{{\bm{T}}}
\newcommand{\bW}{{\bm{W}}}
\newcommand{\bIB}{{\bm{IB}}}
\newcommand{\bOB}{{\bm{OB}}}
\newcommand{\bOS}{{\bm{OS}}}
\newcommand{\bERR}{{\bm{ERR}}}
\newcommand{\bSP}{{\bm{SP}}}
\newcommand{\bMV}{{\bm{MV}}}
\newcommand{\bBM}{{\bm{BM}}}
\newcommand{\balpha}{{\bm{\alpha}}}
\newcommand{\bbeta}{{\bm{\beta}}}
\newcommand{\bgamma}{{\bm{\gamma}}}
\newcommand{\bdelta}{{\bm{\delta}}}
\newcommand{\bomega}{{\bm{\omega}}}
\newcommand{\bzeta}{{\bm{\zeta}}}
\newcommand{\bone}{{\bm{1}}}
\newcommand{\bzero}{{\bm{0}}}

\newcommand{\Cbar}{{\overline{C}}}
\newcommand{\Dbar}{{\overline{D}}}
\newcommand{\dbar}{{\overline{d}}}
\def\Ctilde{{\widetilde{C}}}
\def\Chat{{\widehat{C}}}
\def\ctilde{{\widetilde{c}}}
\def\zbar{{\overline{Z}}}
\def\pitilde{{\widetilde{\pi}}}

%
%
\newcommand{\zfz}{ {{}_0 \! F_0} }
\newcommand{\zfo}{ {{}_0 \! F_1} }
\newcommand{\ofz}{ {{}_1 \! F_0} }
\newcommand{\oft}{ {{}_1 \! F_2} }

%
%
\newcommand{\FHyper}[2]{ {\tensor[_{#1 \!}]{F}{_{#2}}\!} }
\newcommand{\FHYPER}[5]{ {\FHyper{#1}{#2} \!\biggl(
   \!\!\begin{array}{c} #3 \\[1mm] #4 \end{array}\! \bigg|\, #5 \! \biggr)} }
\newcommand{\ofo}{ {\FHyper{1}{1}} }
\newcommand{\tfo}{ {\FHyper{2}{1}} }
\newcommand{\FHYPERbottomzero}[3]{ {\FHyper{#1}{0} \!\biggl(
   \!\!\begin{array}{c} #2 \\[1mm] \hbox{---} \end{array}\! \bigg|\, #3 \! \biggr)} }

%
%
\newcommand{\stirlingsubset}[2]{\genfrac{\{}{\}}{0pt}{}{#1}{#2}}
\newcommand{\stirlingcycleold}[2]{\genfrac{[}{]}{0pt}{}{#1}{#2}}
\newcommand{\stirlingcycle}[2]{\left[\! \stirlingcycleold{#1}{#2} \!\right]}
\newcommand{\assocstirlingsubset}[3]{{\genfrac{\{}{\}}{0pt}{}{#1}{#2}}_{\! \ge #3}}
\newcommand{\genstirlingsubset}[4]{{\genfrac{\{}{\}}{0pt}{}{#1}{#2}}_{\! #3,#4}}
\newcommand{\euler}[2]{\genfrac{\langle}{\rangle}{0pt}{}{#1}{#2}}
\newcommand{\eulergen}[3]{{\genfrac{\langle}{\rangle}{0pt}{}{#1}{#2}}_{\! #3}}
\newcommand{\eulersecond}[2]{\left\langle\!\! \euler{#1}{#2} \!\!\right\rangle}
\newcommand{\eulersecondgen}[3]{{\left\langle\!\! \euler{#1}{#2} \!\!\right\rangle}_{\! #3}}
\newcommand{\binomvert}[2]{\genfrac{\vert}{\vert}{0pt}{}{#1}{#2}}
\newcommand{\binomsquare}[2]{\genfrac{[}{]}{0pt}{}{#1}{#2}}


\newenvironment{sarray}{
             \textfont0=\scriptfont0
             \scriptfont0=\scriptscriptfont0
             \textfont1=\scriptfont1
             \scriptfont1=\scriptscriptfont1
             \textfont2=\scriptfont2
             \scriptfont2=\scriptscriptfont2
             \textfont3=\scriptfont3
             \scriptfont3=\scriptscriptfont3
           \renewcommand{\arraystretch}{0.7}
           \begin{array}{l}}{\end{array}}

\newenvironment{scarray}{
             \textfont0=\scriptfont0
             \scriptfont0=\scriptscriptfont0
             \textfont1=\scriptfont1
             \scriptfont1=\scriptscriptfont1
             \textfont2=\scriptfont2
             \scriptfont2=\scriptscriptfont2
             \textfont3=\scriptfont3
             \scriptfont3=\scriptscriptfont3
           \renewcommand{\arraystretch}{0.7}
           \begin{array}{c}}{\end{array}}

\clearpage

\section{Introduction}

A matrix $M$ of real numbers is called {\em totally nonnegative}\/ (TN)
if every minor of $M$ is nonnegative,
and {\em totally positive}\/ (TP)
if every minor of $M$ is positive.
More generally, $M$ is called {\em totally nonnegative of order $r$}\/
(TN${}_r$)
if every minor of $M$ of size $\le r$ is nonnegative,
and {\em totally positive of order $r$}\/ (TP${}_r$)
if every minor of $M$ of size $\le r$ is positive.\footnote{
   {\bf Warning:}  Some authors
   (e.g.\ \cite{Schoenberg_30,Karlin_68,Ando_87,Pinkus_10,Sokal_totalpos})
   use the terms ``totally positive'' and ``strictly totally positive''
   for what we have termed ``totally nonnegative'' and
   ``totally positive'', respectively.
   So it is very important, when seeing any claim about
   ``totally positive'' matrices, to ascertain which sense of
   ``totally positive'' is being used!
   (This is especially important because many theorems in this subject
    require the {\em strict}\/ concept for their validity:
    see e.g.\ Section~\ref{subsec.TP} below.)
}
Of course, for $m$-by-$n$ matrices,
TN = TN${}_r$ and TP = TP${}_r$ where $r = \min(m,n)$.
Background information on totally nonnegative and totally positive matrices
and their applications can be found in
\cite{Karlin_68,Ando_87,Gasca_96,Gantmacher_02,Pinkus_10,Fallat_11}.

It is an immediate consequence of the Cauchy--Binet formula
that the product of two TN${}_r$ (resp.\ TP${}_r$) matrices
is TN${}_r$ (resp.\ TP${}_r$).
However, other natural matrix operations do not in general
preserve total nonnegativity.
For instance, it is well known (and easy to see by example)
that the sum of two TP matrices need not even be TN${}_2$.
The situation is slightly (but not much) better
when the matrices are symmetric.
Likewise, it has been known for over 40 years that
the Hadamard (entrywise) product of two TN (resp.\ TP) matrices
is always TN${}_2$ (resp.\ TP${}_2$)
but need not be TN${}_3$ \cite[p.~163]{Markham_70}.
Once again, the situation is slightly (but not much) better
when the matrices are symmetric.
In this paper we shall give counterexamples
illustrating the various possibilities
and showing the sharpness of each positive result.

The situation changes radically, however, for Hankel matrices,
i.e.\ square matrices $A = (a_{ij})$ in which $a_{ij}$ depends only on $i+j$.
The Hankel matrices form an important subclass of symmetric matrices,
and they arise in numerous applications
\cite{Shohat_43,Gantmacher_59,Iohvidov_82,Rahman_02,Holtz_12,Sokal_totalpos}.
It is easy to see (Lemma~\ref{lemma.hankel} below) that a matrix is Hankel
if and only if every contiguous submatrix is symmetric.
Here we will exploit this fact to show that,
for Hankel matrices, total nonnegativity ---
and more generally, total nonnegativity of order $r$ ---
{\em is}\/ preserved by sum and by Hadamard product.
We will also show that total nonnegativity of order $r$
is preserved under Hadamard powers with an arbitrary real exponent $t \ge r-2$.

One important motivation for this investigation was the connection between
the Stieltjes moment problem \cite{Shohat_43,Akhiezer_65}
and the total positivity of Hankel matrices.
It is well known that an {\em infinite}\/ Hankel matrix
$A = (a_{i+j})_{i,j=0}^\infty$ is totally nonnegative
if and only if the underlying sequence $\bm{a} = (a_n)_{n=0}^\infty$
is a Stieltjes moment sequence
(i.e.\ the moments of a positive measure on $[0,\infty)$):
``only if'' follows immediately from the standard positive-definiteness
criterion for Stieltjes moment sequences \cite[Theorem~1.3]{Shohat_43},
while ``if'' follows by a simple Vandermonde-matrix argument
\cite[p.~460, Th\'eor\`eme~9]{Gantmakher_37} \cite[Theorem~4.4]{Pinkus_10}.
This equivalence immediately implies that,
for {\em infinite}\/ Hankel matrices,
total nonnegativity is preserved by sum and by Hadamard product.
We therefore wondered whether the same result would hold
when infinite Hankel matrices are replaced by finite ones,
or when TN is replaced by TN${}_r$.
It is satisfying to know that the answer to both questions is yes.

Finally, some of our counterexamples also settle the
Hadamard critical-exponent problem \cite{Johnson_12}
for TN or TP matrices that are general, symmetric or Hankel.

\section{Preliminaries}

In this section we review some known results that will be used
as tools in the remainder of the paper.

\subsection{Inferring total positivity from a proper subset of minors}
   \label{subsec.TP}

We write $[n] = \{1,\ldots,n\}$.
A subset $I \subseteq [n]$ is called {\em contiguous}\/ if it is an interval
(i.e.\ $i,k \in I$ and $i < j < k$ imply $j \in I$).
A subset $I \subseteq [n]$ is called {\em initial}\/ if it is contiguous
and contains 1.

If $A = (a_{ij})_{1 \le i \le m, \, 1 \le j \le n}$ is an $m$-by-$n$ matrix
and $I \subseteq [m]$, $J \subseteq [n]$,
we denote by $A_{IJ}$ the submatrix of $A$
corresponding to the rows $I$ and the columns $J$,
all kept in their original order.
The submatrix $A_{IJ}$ (and the corresponding minor $\det A_{IJ}$)
is called {\em contiguous}\/ if $|I|=|J|$ and both $I$ and $J$ are contiguous;
it is called {\em initial}\/ if $|I|=|J|$ and both $I$ and $J$ are contiguous
and at least one of them is initial.
Note that each matrix entry is the lower-right corner
of exactly one initial submatrix;
so an $m$-by-$n$ matrix has $mn$ initial submatrices.

The following important result \cite[Theorem~4.1]{Gasca_92}
allows one to infer total positivity
from a rather small subset of minors:

\begin{theorem}  {$\!\!\!$ \rm \protect\cite{Gasca_92}\ }
   \label{thm.GP}
Let $A$ be an $m$-by-$n$ matrix.
If all the initial minors of $A$ are positive, then $A$ is TP.
\end{theorem}

\noindent
Proofs can be found in
\cite[Theorem~3.1.4]{Fallat_11} and \cite[Theorem~2.3]{Pinkus_10}.
See also \cite{Fomin_00}
for a combinatorial reinterpretation
of Theorem~\ref{thm.GP} in the case $m=n$,
as well as a generalization to some other sets of minors
of the same cardinality $mn = n^2$.

In fact a weaker result, due to Fekete \cite{Fekete_12} in 1912,
would suffice for our applications:

\begin{theorem}  {$\!\!\!$ \rm \protect\cite{Fekete_12}\ }
   \label{thm.fekete}
Let $A$ be an $m$-by-$n$ matrix.
If all the contiguous minors of $A$ are positive, then $A$ is TP.
\end{theorem}

Please note that Theorems~\ref{thm.GP} and \ref{thm.fekete}
cannot be extended to TN matrices:
for instance \cite[p.~23]{Cryer_76},
all the contiguous minors of the matrix
$A = \displaystyle{ \begin{bmatrix}
                       0 & 1 & 0 \\
                       0 & 0 & 1 \\
                       1 & 0 & 0 \\
                \end{bmatrix}
              }$
are nonnegative, and $\det A = 1$,
but some of the $2$-by-$2$ noncontiguous minors equal~$-1$.
(It~\mbox{follows}~that there is no way of perturbing $A$ so that all
 the initial minors are positive.)
Moreover, the lower-triangular matrix
\cite[p.~86]{Cryer_73} \cite[p.~23]{Cryer_76}
$B = \displaystyle{ \begin{bmatrix}
                       1 & 0 & 0 & 0 \\
                       0 & 1 & 0 & 0 \\
                       0 & 0 & 1 & 0 \\
                       1 & 0 & 0 & 1 \\
                \end{bmatrix}
              }$
and the symmetric matrix
$C = \displaystyle{ \begin{bmatrix}
                       \sqrt{2} & 0 & 0 & 1 \\
                       0 & 1 & 0 & 0 \\
                       0 & 0 & 1 & 0 \\
                       1 & 0 & 0 & \sqrt{2} \\
                \end{bmatrix}
              }$
have the same property.
(Note that the lower-left $3$-by-$3$ corner of $B$ and $C$ is $A$.)
So it does not help here to assume that the matrix is triangular or symmetric.
(But it {\em does}\/ help to assume that the matrix is Hankel:
 see Theorem~\ref{thm.hankel}(a) below.)

\bigskip

In this paper we shall need an extension of
Theorems~\ref{thm.GP} and \ref{thm.fekete}
to total positivity of order $r$,
which we state as follows:

\begin{theorem} {$\!\!\!$ \rm \protect\cite[Corollary~3.1.7]{Fallat_11}\ }
   \label{thm.TPr}
Let $A$ be an $m$-by-$n$ matrix.
Suppose that all the initial minors of~$A$ of size $\le r-1$
are positive, and that
all the contiguous minors of~$A$ of size $r$ are positive.
Then $A$ is TP${}_r$.
\end{theorem}

\noindent
Since this result is stated in \cite{Fallat_11} without proof,
it is perhaps useful to include a proof here.
The first step is to establish the following weakened version
of Theorem~\ref{thm.TPr}:

\begin{lemma}
   \label{lemma.TPr}
Let $A$ be an $m$-by-$n$ matrix.
Suppose that all the contiguous minors of~$A$ of size $\le r$
are positive.
Then $A$ is TP${}_r$.
\end{lemma}

\noindent
The statement above is essentially \cite[Corollary~3.1.6]{Fallat_11},
but for the reader's convenience we give the proof
(which is slightly streamlined compared to the one given in \cite{Fallat_11}):

\proofof{Lemma~\ref{lemma.TPr}}
Let $I \subseteq [m]$ be any contiguous set of rows of cardinality $r$,
and let $J \subseteq [n]$ be any contiguous set of columns of cardinality $r$.
Then Theorem~\ref{thm.fekete} tells us that the
submatrices $A_{I [n]}$ and $A_{[m] J}$ are TP.
Now let $I' = \{i_1,\ldots,i_r\} \subseteq [m]$ be any collection of
$r$ rows, and let $B = A_{I' [n]}$ be the corresponding submatrix.
Note that any $r$-by-$r$ contiguous submatrix of $B$
must lie in some collection of $r$ consecutive columns,
hence be a submatrix of some $A_{[m] J}$.
So any such submatrix of $B$ must be TP;
and applying Theorem~\ref{thm.fekete} again, we deduce that $B$ is TP.
But since any square submatrix of $A$ of size $\le r$
is a submatrix of some such matrix $B$,
we conclude that $A$ is TP${}_r$.
\qed

\proofof{Theorem~\ref{thm.TPr}}
%
By Theorem~\ref{thm.GP},
the submatrix $A_{[r][n]}$ consisting of the initial $r$ rows of $A$ is TP.
Now let $B$ be the submatrix of $A$ obtained by deleting the first row.
Every initial minor of size $\le r-1$ of $B$
is either an initial minor of $A$ or a minor of $A_{[r][n]}$,
hence positive;
and by hypothesis the contiguous minors of size $r$ of $B$
are positive.
So by Theorem~\ref{thm.GP} again,
the submatrix $B_{[r][n]}$ consisting of the initial $r$ rows of $B$ is TP.
Continuing in this manner, we can show that every submatrix of $A$
consisting of $r$ consecutive rows is TP.
Lemma~\ref{lemma.TPr} then implies that $A$ is TP${}_r$.
\qed

%


Our applications will in fact require only Lemma~\ref{lemma.TPr},
not the stronger Theorem~\ref{thm.TPr}.

\subsection{Density of total positivity within total nonnegativity}
   \label{subsec.density}

Since many important properties of TP matrices
(like those in the preceding subsection)
do not extend to TN matrices,
it is very useful to be able to approximate TN matrices by TP matrices,
or more generally to approximate TN${}_r$ matrices by TP${}_r$ matrices.
It is a well-known fact \cite[Theorem~2.6]{Pinkus_10}
that the $m$-by-$n$ TP matrices are dense in the $m$-by-$n$ TN matrices.
The same proof establishes, {\em mutatis mutandis}\/,
the corresponding result for TP${}_r$ and TN${}_r$:

\begin{theorem}
   \label{thm.density}
The set of $m$-by-$n$ TP${}_r$ matrices is dense in the set of
$m$-by-$n$ TN${}_r$ matrices.
\end{theorem}

Since the proof of \cite[Theorem~2.6]{Pinkus_10} preserves symmetry,
we can also assert:

\begin{theorem}
   \label{thm.density.symmetric}
The set of $n$-by-$n$ symmetric TP${}_r$ matrices is dense in the set of
$n$-by-$n$ symmetric TN${}_r$ matrices.
\end{theorem}

The corresponding result for Hankel matrices will be proven in
Corollary~\ref{cor.density.hankel} below.

\subsection{Total positivity of Hankel matrices}
   \label{subsec.hankel}

An $m$-by-$n$ matrix $A = (a_{ij})_{1 \le i \le m, \, 1 \le j \le n}$
is said to be a {\em Hankel matrix}\/ if $a_{ij}$ depends only on $i+j$,
i.e.\ $a_{ij} = a_{i'j'}$ whenever $i+j = i'+j'$.
Hankel matrices are characterized combinatorially
by the following simple but important fact:

\begin{lemma}
   \label{lemma.hankel}
For an $m$-by-$n$ matrix $A$, the following are equivalent:
\begin{itemize}
   \item[(a)] $A$ is Hankel.
   \item[(b)] Every contiguous submatrix of $A$ is Hankel.
   \item[(c)] Every contiguous submatrix of $A$ is symmetric.
   \item[(d)] Every contiguous $2$-by-$2$ submatrix of $A$ is symmetric.
\end{itemize}
\end{lemma}

\proof
(a)$\implies$(b):
Let $A$ be a Hankel matrix,
and consider a contiguous submatrix $A_{IJ}$
where $I = \{i,i+1,\ldots,i+k\}$ and $J = \{j,j+1,\ldots,j+k\}$.
Then the $(s,t)$ entry of $A_{IJ}$ is $a_{i+s-1,j+t-1}$;
and this equals the $(s',t')$ entry whenever $s+t = s'+t'$,
by virtue of the Hankel property of $A$.

(b)$\implies$(c)$\implies$(d) is obvious.

(d)$\implies$(a):
If every contiguous $2$-by-$2$ submatrix of $A$ is symmetric,
we have $a_{i+1,j} = a_{i,j+1}$ for all $i \in [m-1]$ and $j \in [n-1]$.
By combining these facts for different $i,j$,
it is easily seen that $A$ is Hankel.
\qed

In the remainder of this paper, we shall consider only {\em square}\/
Hankel matrices (i.e.\ $m=n$).

In studying Hankel matrices it is often convenient to number
the rows and columns from 0 to $n-1$ rather than from 1 to $n$,
as this facilitates the connection with the Stieltjes moment problem.
Thus, an $n$-by-$n$ Hankel matrix is of the form
$A = (a_{i+j})_{0 \le i,j \le n-1}$
for some sequence of numbers $a_0,a_1,\ldots,a_{2n-2}$.
Assuming that $n \ge 2$,
let us also define $A' = (a_{i+j+1})_{0 \le i,j \le n-2}$,
i.e.\ $A'$ is the $(n-1)$-by-$(n-1)$ submatrix of $A$
in its upper right (or lower left) corner.

The conditions for a finite Hankel matrix to be TN or TP
are slightly delicate, because they involve the theory of
the truncated Stieltjes moment problem \cite{Curto_91}.
But the conditions for an infinite Hankel matrix to be TN or TP
are quite simple, and this is all we shall need here;
indeed, we shall need only the TP case.
Given an infinite sequence $\ba = (a_k)_{k=0}^\infty$ and an integer $m \ge 0$,
let us define the $m$-shifted $n$-by-$n$ Hankel matrix
$H_n^{(m)}(\ba) = (a_{i+j+m})_{0 \le i,j \le n-1}$
and the $m$-shifted infinite Hankel matrix
$H_\infty^{(m)}(\ba) = (a_{i+j+m})_{i,j \ge 0}$.
We then have the following well-known result:

\begin{theorem}
   \label{thm.hankelreal_strict.infty}
For a sequence $\ba = (a_k)_{k=0}^\infty$ of real numbers,
the following are equivalent:
\begin{itemize}
   \item[(a)]  $H^{(0)}_\infty(\ba)$ is totally positive.
      [Equivalently, $H^{(0)}_n(\ba)$ is totally positive for all $n$.]
   \item[(b)]  Both $H^{(0)}_\infty(\ba)$ and $H^{(1)}_\infty(\ba)$
      are positive-definite.
      [Equivalently, $H^{(0)}_n(\ba)$ and $H^{(1)}_n(\ba)$
      are positive-definite for all $n$.]
   \item[(c)]  The leading principal minors
      $\Delta_n = \det H^{(0)}_n(\ba)$ and $\Delta'_n = \det H^{(1)}_n(\ba)$
      are positive for all~$n$.
   \item[(d)]  There exists a positive measure $\mu$ on $[0,\infty)$,
      whose support is an infinite set,
      such that $a_k = \int x^k \, d\mu(x)$ for all $k \ge 0$.
\end{itemize}
\end{theorem}

Indeed, (a)$\implies$(b)$\implies$(c) is trivial,
and (b)$\iff$(c) is Sylvester's criterion for positive-definiteness;
(c)$\iff$(d) is the standard criterion for the Stieltjes moment problem
\cite[Theorem~1.3]{Shohat_43};
and (d)$\implies$(a)
follows by a simple Vandermonde-matrix argument
\cite[p.~460, Th\'eor\`eme~9]{Gantmakher_37} \cite[Theorem~4.4]{Pinkus_10}.

\section{Sums}

The sum of two TN${}_1$ matrices is trivially TN${}_1$,
and the sum of a TN${}_1$ matrix and a TP${}_1$ matrix is trivially TP${}_1$.
But simple examples show that the sum of two TP matrices
need not be TN${}_2$, even if {\em one}\/ of the two matrices is symmetric.

The sum of {\em two}\/ symmetric $2$-by-$2$ TN matrices is TN:
the $1$-by-$1$ minors are covered by the trivial argument,
and the $2$-by-$2$ determinant is nonnegative
because the sum of two positive-semidefinite matrices
is positive-semidefinite.
(Since a $2$-by-$2$ symmetric matrix is automatically Hankel,
this result is a special case of Corollary~\ref{cor.hankel.sum} below.)
But the corresponding assertion fails already for $3$-by-$3$
symmetric matrices, even when one of the two matrices is Hankel:

\begin{example}
   \label{exam.sum3sym}
\rm
\cite[p.~21]{Berman_03}\quad
The symmetric matrix
$I = \displaystyle{ \begin{bmatrix}
                       1 & 0 & 0 \\
                       0 & 1 & 0 \\
                       0 & 0 & 1
                \end{bmatrix}
              }$
and the Hankel matrix
$J = \displaystyle{ \begin{bmatrix}
                       1 & 1 & 1 \\
                       1 & 1 & 1 \\
                       1 & 1 & 1
                \end{bmatrix}
              }$
are TN, but
$I+J = \displaystyle{ \begin{bmatrix}
                       2 & 1 & 1 \\
                       1 & 2 & 1 \\
                       1 & 1 & 2
                \end{bmatrix}
              }$
is not even TN${}_2$.

Moreover, as the reader can easily verify,
these two input matrices can be perturbed slightly to make them TP
while preserving the symmetry and the Hankel property.
\myendremark
\end{example}

But if {\em both}\/ input matrices are Hankel, we have a positive result:

\begin{theorem}
   \label{thm.hankel}
\quad\hfill\vspace*{-1mm}
\begin{itemize}
   \item[(a)]  Let $A$ be a Hankel matrix, all of whose contiguous minors
      of size $\le r$ are nonnegative.  Then $A$ is TN${}_r$.
   \item[(b)]  Let $A$ and $B$ be Hankel matrices,
      all of whose contiguous minors of size $\le r$ are nonnegative.
      Then $A+B$ is TN${}_r$.
   \item[(c)]  Let $A$ (resp.\ $B$) be a Hankel matrix,
      all of whose contiguous minors of size $\le r$ are nonnegative
      (resp.\ positive).
      Then $A+B$ is TP${}_r$.
\end{itemize}
\end{theorem}

\proof
{\bf See Corrigendum (at end of file) for corrected proof of
 Theorem~\ref{thm.hankel}.}

\sout{
Let $A$ be a Hankel matrix, all of whose contiguous minors of size $\le r$
are nonnegative,
and let $B$ be a Hankel matrix, all of whose contiguous minors of size $\le r$ 
are positive.
Then
\mbox{\cite[Theorem~7.2.5]{Horn_13}}
every contiguous submatrix of size $\le r$ of $A$ (resp.~$B$) is
positive-semidefinite (resp.~positive-definite).
It follows that every contiguous submatrix of size $\le r$ of $A+B$
is positive-definite, and hence in particular has a positive
determinant.  By Lemma~\ref{lemma.TPr}, $A+B$ is TP${}_r$.
This proves (c).
}

\sout{
But now, given $A$, we can take $B$ to be any TP Hankel matrix:
that is, by Theorem~\ref{thm.hankelreal_strict.infty}
we can take $B$ to be the Hankel matrix associated to
any Stieltjes moment sequence of infinite support
(for instance, $a_k = k!$ or $a_k = \lambda^{k^2}$ with $\lambda > 1$).
Then $A+\epsilon B$ is TP${}_r$ for all $\epsilon > 0$, hence $A$ is TN${}_r$.
This proves (a).
}

\sout{
Finally, let $A$ and $B$ be Hankel matrices,
all of whose contiguous minors of size $\le r$ are nonnegative.
Then every contiguous submatrix of size $\le r$ of $A$ or $B$ is
positive-semidefinite, so the same holds for $A+B$.
Applying (a) to $A+B$, we obtain (b).
}
\qed

\begin{corollary}
   \label{cor.hankel.sum}
\quad\hfill\vspace*{-1mm}
\begin{itemize}
   \item[(a)]  The sum of two TN${}_r$ Hankel matrices is TN${}_r$.
   \item[(b)]  The sum of a TN${}_r$ Hankel matrix and a TP${}_r$
      Hankel matrix is TP${}_r$.
\end{itemize}
\end{corollary}

\medskip

\begin{corollary}
   \label{cor.density.hankel}
The TP${}_r$ Hankel matrices are dense in the TN${}_r$ Hankel matrices.
\end{corollary}

\medskip

\begin{corollary}
   \label{cor.hankelreal}
Let $A = (a_{i+j})_{0 \le i,j \le n-1}$ be an $n$-by-$n$ Hankel matrix,
and define $A' = (a_{i+j+1})_{0 \le i,j \le n-2}$.
Then $A$ is TN${}_r$ (resp.\ TP${}_r$)
if and only if every principal minor of $A$ and $A'$ of size $\le r$
is nonnegative (resp.\ positive).
In particular, $A$ is TN (resp.\ TP)
if and only if both $A$ and $A'$ are
positive-semidefinite (resp.\ positive-definite).
\end{corollary}

\proofof{Corollary~\ref{cor.hankelreal}}
This is an immediate consequence of Theorem~\ref{thm.hankel}(a)
(resp.\ Lemma~\ref{lemma.TPr}) together with the observation that
every contiguous minor of $A$ is a principal minor of either $A$ or $A'$.
\qed

\bigskip

{\bf Remarks.}
1.  In the same way that Theorem~\ref{thm.TPr}
improves Lemma~\ref{lemma.TPr},
one might hope to improve Theorem~\ref{thm.hankel}
by weakening the hypothesis on contiguous minors of size $\le r$ to
``initial minors of size $\le r-1$ and contiguous minors of size $r$''.
But it turns out that this does {\em not}\/ hold in general.
Consider, for any $n \ge 3$, the $n$-by-$n$ matrix $A$
having $1$ in the upper-left corner, $-1$ in the lower-right corner,
and zeros elsewhere.
Then $A$ is Hankel;
all its initial minors of size 1 are either 0 or 1,
all its contiguous minors of size 2 are zero,
and all its minors of size $\ge 3$ are zero;
but $A$ is not even TN${}_1$, much less TN${}_r$ for some or all $r \in [2,n]$.

2.  In a general partially ordered commutative ring,
the sum of two TN Hankel matrices can {\em fail}\/ to be TN,
even if one of the two matrices is a matrix of pure numbers:
for instance, in the polynomial ring $\R[x]$ with the coefficientwise order,
$A = \displaystyle{ \begin{bmatrix}
                       1 & 1 \\
                       1 & 1
                \end{bmatrix}
              }$
and
$B = \displaystyle{ \begin{bmatrix}
                       1 & x \\
                       x & x^2
                \end{bmatrix}
              }$
are TN,
but $\det(A+B) = 1 - 2x + x^2$ fails to be coefficientwise nonnegative.
See \cite{Sokal_totalpos} for further discussion.
\myendremark

\bigskip

Finally, let us return to general matrices,
and pose the following question:
Which $m$-by-$n$ matrices $A$ have the property that
$A+B$ is TN whenever $B$ is TN?
The answer is as follows:

\begin{theorem}
   \label{thm.addcore}
Let $A = (a_{ij})$ be an $m$-by-$n$ matrix.
Then the following are equivalent:
\begin{itemize}
   \item[(a)]  $A+B$ is TN whenever $B$ is TN.
   \item[(b)]  $A+B$ is TP whenever $B$ is TP.
   \item[(c${}_r$)]  $A+B$ is TN${}_r$ whenever $B$ is TN${}_r$.
     [Here $r \ge 2$.]
   \item[(d${}_r$)]  $A+B$ is TP${}_r$ whenever $B$ is TP${}_r$.
     [Here $r \ge 2$.]
   \item[(e)]  $A+B$ is TN${}_2$ whenever $B$ is TN.
   \item[(f)]  $A+B$ is TN${}_2$ whenever $B$ is TP.
   \item[(g)]  The upper-left and lower-right elements of $A$ are nonnegative,
      and all other elements are zero.
\end{itemize}
\end{theorem}

\proof
Suppose first that all the elements of $A$ are zero except possibly
$a_{11}$ and $a_{mn}$,
and consider an $r$-by-$r$ minor $\det\, (A+B)_{IJ}$
according to whether it contains the first row and column,
the last row and column, both, or neither.
If neither, then obviously $\det\, (A+B)_{IJ} = \det B_{IJ}$.
If the first but not the last,
then $\det\, (A+B)_{IJ} = \det B_{IJ} +
  a_{11} \det B_{I \setminus 1, J \setminus 1}$.
If the last but not the first,
then $\det\, (A+B)_{IJ} = \det B_{IJ} +
  a_{mn} \det B_{I \setminus m, J \setminus n}$.
And if both, then
then $\det\, (A+B)_{IJ} = \det B_{IJ} +
  a_{11} \det B_{I \setminus 1, J \setminus 1} +
  a_{mn} \det B_{I \setminus m, J \setminus n} +
  a_{11} a_{mn} \det B_{I \setminus \{1,m\}, J \setminus \{1,n\}}$.
It follows that (g) implies all the other statements.

Conversely, suppose that $A+B$ is TN${}_2$ whenever $B$ is TN.
Taking $B=0$, we conclude that $A$ must be TN${}_2$
and in particular all its entries must be nonnegative.
Furthermore, if $(i,j)$ is any entry other than $(1,1)$ or $(m,n)$,
then we can choose $(i',j')$ to be either $(i-1,j+1)$ or $(i+1,j-1)$
and let $B$ be the matrix with $b_{i'j'} = \lambda$ and all other entries zero.
Then $B$ is TN whenever $\lambda \ge 0$, and we have
$\det (A+B)_{\{i,i'\},\{j,j'\}} =
 \det A_{\{i,i'\},\{j,j'\}} - \lambda a_{ij}$.
Taking $\lambda \to +\infty$ we conclude that $a_{ij} \le 0$.
This proves (e) $\implies$ (g).
Finally, (f) $\implies$ (e)
is an immediate consequence of Theorem~\ref{thm.density}.
\qed

The matrices $A$ characterized in Theorem~\ref{thm.addcore}
could be termed the ``additive core'' of the TN matrices,
by analogy with the ``Hadamard core'' studied in \cite{Crans_01}.

\section{Hadamard product}  \label{sec.hadamard}

If $A = (a_{ij})$ and $B = (b_{ij})$ are two matrices of the same size
(say, $m$-by-$n$),
their {\em Hadamard product}\/ (or {\em entrywise product}\/) $A \circ B$
is the matrix with elements $(A \circ B)_{ij} = a_{ij} b_{ij}$.
See \cite[Chapter~5]{Horn_91} for further information on the properties
of the Hadamard product.

The Hadamard product of two TN${}_1$ (resp.\ TP${}_1$) matrices
is trivially TN${}_1$ (resp.\ TP${}_1$).
Moreover, the TN${}_2$ and TP${}_2$ cases
are handled by the following easy positive result:

\begin{lemma}
   \label{lemma.hadamard_product}
\quad\hfill\vspace*{-1mm}
\begin{itemize}
   \item[(a)] The Hadamard product of two TN${}_2$ matrices is TN${}_2$.
   \item[(b)] The Hadamard product of a TP${}_2$ matrix with a
      TP${}_1$ $\cap$ TN${}_2$ matrix is TP${}_2$.
\end{itemize}
\end{lemma}

%
%

However, in general the Hadamard product of two TP matrices
need not even be TN${}_3$,
as was observed
in \cite[p.~163]{Markham_70}.
The precise situation is as follows.

First, an easy positive result concerning $3$-by-$3$ {\em symmetric}\/
matrices:

\begin{proposition}
   \label{prop.hadamard3sym}
The Hadamard product of two $3$-by-$3$ TN (resp.\ TP) symmetric matrices
is TN (resp.\ TP).
\end{proposition}

\proof
The $1$-by-$1$ minors are trivially nonnegative (resp.\ positive).
The $2$-by-$2$ minors are nonnegative (resp.\ positive)
by Lemma~\ref{lemma.hadamard_product}.
The $3$-by-$3$ determinant is nonnegative (resp.\ positive)
by the Schur product theorem \cite[Theorem~7.5.3]{Horn_13}.
\qed

But the corresponding result fails for $3$-by-$3$ nonsymmetric matrices,
and for $4$-by-$4$ symmetric matrices, as the following two examples
demonstrate:

\begin{example}
   \label{exam.hadamard3asym}
\rm
\cite{Crans_01}\quad
The $3$-by-$3$ matrices
$W = \displaystyle{ \begin{bmatrix}
                        1 & 1 & 0 \\
                        1 & 1 & 1 \\
                        1 & 1 & 1
                    \end{bmatrix}
                  }$
and
$W^{\rm T} = \displaystyle{ \begin{bmatrix}
                        1 & 1 & 1 \\
                        1 & 1 & 1 \\
                        0 & 1 & 1
                    \end{bmatrix}
                  }$
are TN, but
$W \circ W^{\rm T} = \displaystyle{ \begin{bmatrix}
                        1 & 1 & 0 \\
                        1 & 1 & 1 \\
                        0 & 1 & 1
                    \end{bmatrix}
                  }$
is not (its determinant is~$-1$).
Moreover, by a suitable small perturbation we can take
the two starting matrices to be TP.

It also does not help to assume that {\em one}\/ of the two matrices
is symmetric or even Hankel, when the other is nonsymmetric:
for instance, if
$A = \displaystyle{ \begin{bmatrix}
                        a_0 & 2 & 1 \\
                        2 &  1 & 1 \\
                        1 &  1 & 2
                    \end{bmatrix}
                  }$
and $W$ is as above,
then $A$ is TP whenever $a_0 > 5$,
but $A \circ W$ is TN only when $a_0 \ge 6$.

Moreover, by ``exterior bordering'' \cite[Theorem~9.0.1]{Fallat_11}
there exist {\em arbitrarily large}\/ TP matrices $A$
since that $A \circ A^{\rm T}$ is not even TN${}_3$.
\myendremark
\end{example}

\begin{example}
   \label{exam.hadamard4sym}
\rm
The $4$-by-$4$ symmetric matrices
\be
   A \;=\; \begin{bmatrix}
               2 & 2 & 1 & 1 \\
               2 & 2 & 1 & 1 \\
               1 & 1 & 2 & 2 \\
               1 & 1 & 2 & 2
           \end{bmatrix}
   \quad\hbox{and}\quad
   B \;=\; \begin{bmatrix}
              2 & 1 & 1 & 0 \\
              1 & 2 & 2 & 1 \\
              1 & 2 & 2 & 1 \\
              0 & 1 & 1 & 2
           \end{bmatrix}
\ee
are TN, but their Hadamard product
\be
   A \circ B
   \;=\;
   \begin{bmatrix}
       4 & 2 & 1 & 0 \\
       2 & 4 & 2 & 1 \\
       1 & 2 & 4 & 2 \\
       0 & 1 & 2 & 4
   \end{bmatrix}
\ee
has several negative $3$-by-$3$ minors
(for instance, the upper right $3$-by-$3$ submatrix
 has determinant $-6$).
Moreover, by a suitable small perturbation
(using Theorem~\ref{thm.density.symmetric})
we can take the two starting matrices to be TP.

It also does not help to assume that {\em one}\/ of the two symmetric
matrices is Hankel, when the other is not.
For instance, the Hankel matrix
\be
   H
   \;=\;
   \begin{bmatrix}
       a_0 & a_1 & 2 & 1 \\
       a_1 & 2   & 1 & 2 \\
       2   & 1   & 2 & a_1 \\
       1   & 2   & a_1 & a_0
   \end{bmatrix}
\ee
is TP whenever $a_1 > 7$ and $a_0 > a_1^2 - 4a_1 + 5$.
But
\be
   A \circ H
   \;=\;
   \begin{bmatrix}
       2a_0 & 2a_1 & 2 & 1 \\
       2a_1 & 4   & 1 & 2 \\
       2   & 1   & 4 & 2a_1 \\
       1   & 2   & 2a_1 & 2a_0
   \end{bmatrix}
\ee
fails to be TN${}_3$ whenever
$a_1 < 4 + \frac{3}{2} \sqrt{5} \approx 7.354102$;
and similarly $B \circ H$ fails to be TN${}_3$ whenever
$a_1 < 4 + 2\sqrt{3} \approx 7.464102$.
\myendremark
\end{example}

But if {\em both}\/ input matrices are Hankel, we again have a positive result:

\begin{theorem}
   \label{thm.hankel.hadamard_product}
\quad\hfill\vspace*{-1mm}
\begin{itemize}
   \item[(a)]  The Hadamard product of two TN${}_r$ Hankel matrices is TN${}_r$.
   \item[(b)]  The Hadamard product of two TP${}_r$ Hankel matrices is TP${}_r$.
\end{itemize}
\end{theorem}


\proof
Suppose that $A$ and $B$ are TN${}_r$ (resp.\ TP${}_r$) Hankel matrices.
Then
all their contiguous submatrices of size $\le r$ are
positive-semidefinite (resp.\ positive-definite).
Therefore, by the Schur product theorem \cite[Theorem~7.5.3]{Horn_13},
all the contiguous submatrices of $A \circ B$ of size $\le r$
are positive-semidefinite (resp.\ positive-definite)
and in particular have a nonnegative (resp.\ positive) determinant.
By Theorem~\ref{thm.hankel}(a) (resp.\ Lemma~\ref{lemma.TPr})
it follows that $A \circ B$ is TN${}_r$ (resp.\ TP${}_r$).
[Alternatively, one could first prove (b) and then invoke
 Corollary~\ref{cor.density.hankel} to deduce (a).]
\qed

{\bf Remarks.}
1.  In a general partially ordered commutative ring
--- for instance, in the polynomial ring $\R[x]$ with the
coefficientwise order ---
the Hadamard product of two TN Hankel matrices can {\em fail}\/ to be TN,
even if one of the two matrices is a matrix of pure numbers.
Furthermore, the Hadamard square of a TN Hankel matrix can fail to be TN.
See \cite{Sokal_totalpos} for details.

2. Many further results concerning total positivity and the Hadamard product
can be found in \cite[Chapter~8]{Fallat_11}.
For instance, the Hadamard product of a TN matrix
and a {\em tridiagonal}\/ TN matrix is TN \cite[Theorem~8.2.5]{Fallat_11};
and this result extends to TN${}_r$, by the same proof.
\myendremark

\section{Hadamard powers}

If $A = (a_{ij})$ is a matrix
and $t > 0$ is an integer,
the {\em Hadamard power}\/ (or {\em entrywise power}\/) $A^{\circ t}$
is defined to be the matrix with elements $(A^{\circ t})_{ij} = a_{ij}^t$.
Moreover, if the matrix $A$ has {\em nonnegative}\/ real entries
--- as we shall assume henceforth ---
we can make this same definition for arbitrary {\em real}\/ powers $t > 0$.

Note that each minor of $A^{\circ t}$ is an exponential polynomial
$ f(t)  =  \sum_{i=1}^n a_i \, e^{\lambda_i t} $,
where we can assume that $a_1,\ldots,a_n$ are real and nonzero
and $\lambda_1 < \ldots < \lambda_n$.
Laguerre's rule of signs
\cite{Laguerre_1883} \cite[pp.~46--47, Problem~V.77]{Polya-Szego}
\cite{Jameson_06}
then states that the number of real zeros of $f$ (counting multiplicity)
is at most the number of sign changes
in the sequence $a_1,\ldots,a_n$,
and is also of the same parity.

If $A$ is TN${}_1$ (resp.\ TP${}_1$),
then trivially so is $A^{\circ t}$ for all real $t > 0$.
Moreover, the following result is almost trivial:

\begin{proposition}
   \label{prop.TN2.hadpow}
If $A$ is TN${}_2$ (resp.\ TP${}_2$),
then so is $A^{\circ t}$ for all real $t > 0$.
\end{proposition}


Less trivially, the TN${}_3$/TP${}_3$ case
is handled by the following result \cite[Theorem~4.2]{Johnson_12}
(see also \cite[pp.~179--180]{Fallat_11}):
%
%

\begin{theorem}
   \label{thm.JW}
If $A$ is TN${}_3$ (resp.\ TP${}_3$),
then so is $A^{\circ t}$ for all real $t \ge 1$.
%
\end{theorem}

\noindent
Here is a slightly simplified proof:

\smallskip

\proof
It obviously suffices to prove the result for $3$-by-$3$ matrices;
and it suffices to prove the TP${}_3$ case, since the TN${}_3$ case
then follows by Theorem~\ref{thm.density}.
By row and column rescalings it suffices to consider
$A = \displaystyle{ \begin{bmatrix}
                       1 & 1 & 1 \\
                       1 & a & b \\
                       1 & c & d
                   \end{bmatrix}
                 }$.
Such a matrix is TP${}_2$ if and only if $a > 1$, $b > a$, $c > a$
and $ad > bc$.
Then $\det(A^{\circ t}) = (ad)^t - d^t - (bc)^t + b^t + c^t - a^t$,
which equals $(\log a \log d - \log b \log c)t^2 + O(t^3)$ near $t=0$
and hence has at least a double root there;
moreover, $\det(A^{\circ t}) \to +\infty$ as $t \to +\infty$.
By Laguerre's rule of signs,
$f(t) = \det(A^{\circ t})$ has precisely three real roots
(note that the unknown ordering of $b$ and $c$ plays no role here
 because their coefficients in $f$ have the same sign;
 likewise for $d$ and $bc$).
And regardless of whether the third root lies at $t < 0$, $t=0$ or $t > 0$
(which depends on the sign of $\log a \log d - \log b \log c$),
$f(1) > 0$ implies $f(t) > 0$ for all $t > 1$.
\qed

The following example shows that Theorem~\ref{thm.JW}
does not extend to $0 < t < 1$, even if the matrix is assumed to be Hankel:

\begin{example}
   \label{exam.JW}
\rm
The matrix
$A = \displaystyle{ \begin{bmatrix}
                        2 & 3 & 5 \\
                        3 & 5 & 9 \\
                        5 & 9 & 17
                    \end{bmatrix}
                  }$
is the $3$-by-$3$ Hankel matrix associated to the
Stieltjes moment sequence $a_n = 1^n + 2^n$ ($n \ge 0$):
all the $1$-by-$1$ and $2$-by-$2$ minors are positive,
and $\det A = 0$.  But
\be
   \det(A^{\circ t})
   \;=\;
   170^t - 162^t - 153^t + 2 \cdot 135^t - 125^t
   \;,
\ee
which is strictly negative for $0 < t < 1$.\footnote{
   {\sc Proof:}
   For $f(t) = \det(A^{\circ t})$,
   easy computations show that
   $t=0$ is a double root and $t=1$ is a simple root.
   Since the coefficient sequence has three sign changes,
   Laguerre's rule of signs
   implies that $f$ has no other real zeros.
   Moreover, straightforward computations show that $f(t) < 0$
   for $t$ slightly less than 1 and for $t$ slightly greater than 0;
   therefore $f(t) < 0$ for all $t \in (0,1)$.
}

Moreover, by a small perturbation
(using Corollary~\ref{cor.density.hankel})
we can make the matrix Hankel and TP
and have $\det(A^{\circ t}) < 0$ for $\delta < t < 1-\delta$,
for arbitrarily small $\delta > 0$.
\myendremark
\end{example}

%


Let us now present some positive results for symmetric and Hankel matrices.
Our key tool will be the following \cite[Theorem~2.2]{Fitzgerald_77}:

\begin{theorem} {$\!\!\!$ \rm \protect\cite{Fitzgerald_77}\ }
   \label{thm.FH}
Let $n \ge 2$, and let $A$ be a symmetric
positive-semidefinite (resp.\ positive-definite)
$n$-by-$n$ matrix with nonnegative real entries.
Then, for all real $t \ge n-2$ ($t > 0$ if $n=2$),
the Hadamard power $A^{\circ t}$ is
positive-semidefinite (resp.\ positive-definite).
\end{theorem}

\noindent
The non-strict (``positive-semidefinite'') version of this result
was proven in \cite{Fitzgerald_77}.\footnote{
   {\bf Warning:}  In \cite{Fitzgerald_77}
   the term ``positive definite'' is used for what is ordinarily called
   ``positive-semidefinite''.
}
By a perturbation argument one can then deduce the strict version.\footnote{
   Given a positive-definite matrix $A$, choose $\epsilon > 0$
   such that $A - \epsilon I$ is positive-semidefinite.
   Then $A^{\circ t} = (A - \epsilon I)^{\circ t} + D$
   where $D = \diag \bigl(a_{ii}^t - (a_{ii}-\epsilon)^t \bigr)$
   is positive-definite.
}

Let us observe that the bound $t \ge n-2$ in Theorem~\ref{thm.FH}
cannot be improved, even if the matrix is Hankel and TP,
as the following example shows:

\begin{example}
   \label{exam.FH}
\rm
Let $n \ge 2$ and $u_1,\ldots,u_n \in \R$,
and define for each $\epsilon > 0$
the $n$-by-$n$ matrix $A_n(\epsilon) = (1+ \epsilon u_i u_j)_{i,j=1}^n$.
This matrix is symmetric and positive-semidefinite;
if the $u_1,\ldots,u_n$ are not all equal, it is of rank 2
(otherwise it is of rank 1);
and if $0 \le u_1 < u_2 < \ldots < u_n$,
then it is also TN $\cap$ TP${}_2$.
A straightforward computation (see Appendix~\ref{app.B}) shows that
\be
   \det A_n(\epsilon)^{\circ t}
   \;=\;
   \biggl( \prod_{k=1}^{n-1} {1 \over k!} \biggr)
   \biggl( \prod_{1 \le i < j \le n} \! (u_i - u_j)^2 \biggr)
   \biggl( \prod_{k=0}^{n-2} (t-k)^{n-1-k} \biggr)
   \epsilon^{n(n-1)/2}
   \,+\,
   O(\epsilon^{n(n-1)/2 \,+\, 1})
   \;.
 \label{eq.exam.FH}
\ee
Therefore, if $n \ge 3$ and the $u_1,\ldots,u_n$ are all distinct,
then for any $t \in (n-3, n-2)$ we have
$\det A_n(\epsilon)^{\circ t} < 0$ for all sufficiently small $\epsilon > 0$.
More generally, if $t \in (m-3, m-2)$ for some integer $m \in [3,n]$,
then $A_m(\epsilon)^{\circ t}$
--- which is a leading principal submatrix of $A_n(\epsilon)^{\circ t}$ ---
has a negative determinant for small $\epsilon > 0$.
So, for all noninteger $t \in (0,n-2)$,
$A_n(\epsilon)^{\circ t}$ fails to be positive-semidefinite
for small $\epsilon > 0$ (how small may depend on $t$).

In \cite[p.~636]{Fitzgerald_77} the authors chose $u_i = i$
and proved the failure of positive-semidefinite\-ness for
noninteger $t \in (0,n-2)$ and small $\epsilon$ by a different method
(computing the inner product
${\bf x}^{\rm T} A_n(\epsilon)^{\circ t} \, {\bf x}$
in power series in $\epsilon$
for a suitably chosen vector ${\bf x} \in \R^n$).

On the other hand, if we choose $u_i = 2^{i-1}$,
then the matrix $A_n(\epsilon)$ is Hankel and TN $\cap$ TP${}_2$.
Moreover, by a small perturbation we can make the matrix Hankel and TP
(by Corollary~\ref{cor.density.hankel}).
Therefore, for each $n \ge 3$ and each noninteger $t \in (0,n-2)$,
there exists an $n$-by-$n$ TP Hankel matrix $A$ such that
one of the leading principal minors of $A^{\circ t}$ is negative,
and in particular $A^{\circ t}$ fails to be positive-semidefinite.
\myendremark
\end{example}

Returning now to TN and TP for Hadamard powers,
we have the following positive result
for $4$-by-$4$ {\em symmetric}\/ matrices:

\begin{proposition}
   \label{prop.hadprod4sym}
If $A$ is a $4$-by-$4$ TN (resp.\ TP) symmetric matrix,
then so is $A^{\circ t}$ for all real $t \ge 2$.
%
\end{proposition}

\proof
The $1$-by-$1$ minors are handled by the trivial argument;
the $2$-by-$2$ minors are handled by Proposition~\ref{prop.TN2.hadpow};
the $3$-by-$3$ minors are handled by Theorem~\ref{thm.JW};
and the $4$-by-$4$ determinant is handled by Theorem~\ref{thm.FH}.
\qed

Example~\ref{exam.FH} shows that Proposition~\ref{prop.hadprod4sym}
cannot be extended to any $t \in (0,1) \cup (1,2)$,
even if the matrix is Hankel and TP.
It follows, in the language of \cite{Johnson_12},
that the Hadamard critical exponent for $4$-by-$4$ {\em symmetric}\/
TN or TP matrices is 2.

\bigskip

But for $4$-by-$4$ nonsymmetric matrices,
and for $5$-by-$5$ symmetric matrices,
even the Hadamard square ($t=2$) does not in general preserve
total nonnegativity, as we now proceed to show.

In \cite[Example~1]{Fallat_07}
an example was given of a $4$-by-$4$ TP matrix
whose Hadamard square has a negative determinant:
\be
   A  \;=\;
   \begin{bmatrix}
       1 & 11 & 22 & 20 \\
       6 & 67 & 139 & 140 \\
       126 & 182 & 395 & 445 \\
       12 & 138 & 309 & 376
   \end{bmatrix}
   \;,
\ee
which has $\det(A \circ A) = -114904113$.
Here is another example of the same phenomenon,
in which moreover $\det(A^{\circ t}) < 0$ also for real $t > 1$:

%
%
%

\begin{example}
   \label{exam.hadamard4asym}
\rm
Consider the matrix
\be
   A  \;=\;
   \begin{bmatrix}
        1 & 1 & 1 & 1 \\
        1 & 1+x & 1+2x & 1+3x \\
        1 & 1+2x & 1+4x & 1+6x \\
        1 & 1+3x & 1+8x & 1+14x
   \end{bmatrix}
   \;,
\ee
which fails to be symmetric only because $a_{34} \neq a_{43}$.
All the $2$-by-$2$ minors are of the form $a x+ b x^2$
with $a > 0$ and $b \ge 0$;
all the $3$-by-$3$ minors are of the form $c x^2$ with $c \ge 0$;
and $\det(A) = 0$.
So $A$ is coefficientwise TN in the polynomial ring $\R[x]$;
in~particular, it is TN for all $x \ge 0$.
But $\det(A \circ A) = -16x^4 + 248x^5$,
so $\det(A \circ A) < 0$ whenever $0 < x < 2/31$.
Furthermore, for small $x$ we have
$\det(A^{\circ t}) = 2(t^3-t^4) x^4 + O(x^5)$;
so for every real $t > 1$ there exists $\delta_t > 0$
such that $\det(A^{\circ t}) < 0$ whenever $0 < x < \delta_t$.
\myendremark
\end{example}

And by perturbing a few coefficients, TN can be upgraded to TP:

\begin{example}
   \label{exam.hadamard4asym.TP}
\rm
Consider the matrix
\be
   A  \;=\;
   \begin{bmatrix}
        1 & 1 & 1 & 1 \\
        1 & 1+x & 1+2x & 1+3x \\
        1 & 1+2x & 1+ (4+\epsilon)x & 1+ (6+{5 \over 2}\epsilon)x \\
        1 & 1+3x & 1+8x & 1+(14+\epsilon)x
   \end{bmatrix}
\ee
with $0 < \epsilon < 1$.
All the $2$-by-$2$ minors are again of the form $a x+b x^2$
with $a > 0$ and $b \ge 0$;
but now all the $3$-by-$3$ minors are of the form $c x^2$ with $c > 0$,
and $\det A = \epsilon^2 x^3$.
So $A$ is TP for all $x > 0$ and $\epsilon \in (0,1)$.
And for small $x$,
\be
   \det(A^{\circ t})  \;=\;
   \epsilon^2 t^3 x^3
   \:+\:
   {1 \over 4} (8  - 70 \epsilon  - 59 \epsilon^2  - 4 \epsilon^3)
       (t^3 - t^4) x^4
   \:+\:  O(x^5)
   \;.
\ee
Therefore, for every real $t > 1$
and small enough $\epsilon > 0$ (depending on $t$)
there exists a nonempty interval of $x > 0$ such that $\det(A^{\circ t}) < 0$.
\myendremark
\end{example}

Examples~\ref{exam.hadamard4asym} and \ref{exam.hadamard4asym.TP}
answer an open question from \cite[p.~81]{Johnson_12},
by showing that the Hadamard critical exponent for
$4$-by-$4$ TN or TP matrices is $\infty$.
Moreover, by ``exterior bordering'' \cite[Theorem~9.0.1]{Fallat_11}
there exists, for any $n \ge 4$ and $t > 1$,
a TP $n$-by-$n$ matrix $A$ such that $A^{\circ t}$
is not even TN${}_4$.
Therefore, for all $n \ge 4$ and $r \ge 4$,
the Hadamard critical exponent for
$n$-by-$n$ TN${}_r$ or TP${}_r$ matrices is also $\infty$.

And here is an example that is almost Hankel:

%
%
%

\begin{example}
   \label{exam.hadamard4asym_bis}
\rm
Consider the matrix
\be
   A  \;=\;
   \begin{bmatrix}
        1+3x  & 1+6x   & 1+14x  & 1+36x \\
        1+6x  & 1+14x  & 1+36x  & 1+98x \\
        1+14x & 1+36x  & 1+98x  & 1+276x \\
        1+36x & 1+98x  & 1+284x & 1+842x
   \end{bmatrix}
   \;.
\ee
All the $2$-by-$2$ minors are of the form $a x+ b x^2$ with $a, b > 0$;
all the $3$-by-$3$ minors are of the form $c x^2 + d x^3$ with $c,d > 0$;
and $\det(A) = 0$.
So $A$ is coefficientwise TN in the polynomial ring $\R[x]$
and is TN for all $x \ge 0$.
But for small $x$ we have
$\det(A^{\circ t}) = 28584 (t^3-t^4) x^4 + O(x^5)$;
so for every real $t > 1$ there exists $\delta_t > 0$
such that $\det(A^{\circ t}) < 0$ whenever $0 < x < \delta_t$.

Note how this example is constructed:
we start from the $4$-by-$4$ Hankel matrix
$(1 + \alpha_{i+j}x)_{0 \le i,j \le 3}$
associated to the Stieltjes moment sequence $\alpha_n = 1^n + 2^n + 3^n$,
which is coefficientwise TN in $\R[x]$
by a general result in \cite{Sokal_totalpos};
we then modify this matrix by changing $a_{32}$ and $a_{33}$.

By replacing 842 by $842+\epsilon$ in the lower-right matrix entry,
TN can be upgraded to TP analogously to Example~\ref{exam.hadamard4asym.TP}.
\myendremark
\end{example}

We can now exhibit a $5$-by-$5$ symmetric TP matrix
whose Hadamard powers with $t > 1$ fail to be TN${}_4$;
indeed, we will choose this $5$-by-$5$ symmetric matrix
so that the $4$-by-$4$ submatrix in its upper right corner
is precisely the almost-Hankel matrix of Example~\ref{exam.hadamard4asym_bis}:

\begin{example}
   \label{exam.hadamard5sym}
\rm
Consider the matrix
\be
   A  \;=\;
   \begin{bmatrix}
        1+2x  & 1+3x  & 1+6x   & 1+14x  & 1+36x \\
        1+3x  & 1+6x  & 1+14x  & 1+36x  & 1+98x \\
        1+6x  & 1+14x & 1+36x  & 1+98x  & 1+276x \\
        1+14x & 1+36x & 1+98x  & 1+284x & 1+842x \\
        1+36x & 1+98x & 1+276x & 1+842x & 1+2604x
   \end{bmatrix}
   \;.
\ee
All the $2$-by-$2$ minors are of the form $ax+bx^2$ with $a,b > 0$;
all the $3$-by-$3$ minors are of the form $cx^2 + dx^3$ with $c,d > 0$;
all the $4$-by-$4$ minors are of the form $ex^3 + fx^4$ with $e,f \ge 0$;
and $\det(A) = 0$.
So $A$ is coefficientwise TN in $\R[x]$ and is TN for all $x \ge 0$.
But the $4$-by-$4$ upper-right submatrix of $A^{\circ t}$
has a negative determinant in the circumstances
discussed in Example~\ref{exam.hadamard4asym_bis}.
(Once again, TN can be upgraded to TP by a small perturbation.)
\myendremark
\end{example}

Example~\ref{exam.hadamard5sym} shows that
the Hadamard critical exponent for
$5$-by-$5$ symmetric TN or TP matrices is $\infty$;
and by ``exterior bordering'' the same result holds for
$n$-by-$n$ symmetric TN${}_r$ or TP${}_r$ matrices
for all $n \ge 5$ and $r \ge 4$.

Finally, for Hankel matrices we have the following positive result:

\begin{theorem}
   \label{thm.hankel.hadamard_power}
For every integer $r \ge 3$:
If $A$ is a TN${}_r$ (resp.\ TP${}_r$) Hankel matrix,
then so is $A^{\circ t}$ for all real $t \ge r-2$.
%
\end{theorem}

\proof
By Theorem~\ref{thm.hankel}(a) (resp.\ Lemma~\ref{lemma.TPr}),
it suffices to show that all the contiguous minors of $A^{\circ t}$
of size $\le r$ are nonnegative (resp.\ positive).
Since all the contiguous submatrices of $A$ of size $\le r$
are symmetric and positive-semidefinite (resp.\ positive-definite),
this is an immediate consequence of Theorem~\ref{thm.FH}.
\qed

Example~\ref{exam.FH} shows that for every $n \ge r$
and every noninteger $t \in (0,r-2)$,
there is an $n$-by-$n$ TP Hankel matrix $A$ such that
one of the leading principal minors of $A^{\circ t}$ of size $\le r$
is negative.
So the bound $t \ge r-2$ in Theorem~\ref{thm.hankel.hadamard_power}
cannot be improved.
In other words, the Hadamard critical exponent
for TN${}_r$ or TP${}_r$ Hankel matrices is exactly $r-2$.

\bigskip

{\bf Final remark.}
Most of the counterexamples in this paper were found by applying
{\sc Mathematica}'s function {\tt FindInstance} to a suitably chosen Ansatz,
sometimes followed by experimentation to find a simpler ``nearby'' example.
Since in practice this works (with present-day hardware and software)
only if the Ansatz has at most three or four parameters,
considerable trial and error was sometimes needed
to find a suitable Ansatz.

\appendix

\section{Proof of equation~\reff{eq.exam.FH}} \label{app.B}

Let $n \ge 2$ and $u_1,\ldots,u_n \in \C$,
and set $M = \max\limits_{1 \le i \le n} |u_i|$.
We consider the $n$-by-$n$ matrix
$A_n(\epsilon) = (1+ \epsilon u_i u_j)_{i,j=1}^n$ for $\epsilon \in \C$.
The binomial series for $(1 + \epsilon u_i u_j)^t$ is convergent
for $|\epsilon| < 1/M^2$ and yields
\be
   A_n(\epsilon)^{\circ t}
   \;=\;
   \sum_{k=0}^\infty \epsilon^k \, {t^{\underline{k}} \over k!} \,
         {\bf u}^{[k]} {\bf u}^{[k] \rm T}
\ee
where $t^{\underline{k}} = t (t-1) \cdots (t-k+1)$
and ${\bf u}^{[k]} = (u_1^k, \ldots, u_n^k)$.
We therefore have $A_n(\epsilon)^{\circ t} = V D V^{\rm T}$
where $V = (u_i^j)_{1 \le i \le n, \, j \ge 0}$
is the $n$-by-$\infty$ Vandermonde matrix,
and $D = \diag(\epsilon^j t^{\underline{j}} /j!)_{j \ge 0}$
is a diagonal matrix.
The Cauchy--Binet formula then gives
\be
   \det A_n(\epsilon)^{\circ t}
   \;=\;
   \sum_J (\det V_{[n] J})^2 \,
       \prod_{j \in J} {\epsilon^j t^{\underline{j}} \over j!}
   \;,
\ee
where the sum runs over $n$-element subsets $J \subseteq \N$.
The smallest power of $\epsilon$ comes from $J = \{0,1,\ldots,n-1\}$
and in this case $\det V_{[n] J}$ is the Vandermonde determinant
$\prod\limits_{1 \le i < j \le n} (u_j - u_i)$.
All other terms contribute higher powers of $\epsilon$
(and the coefficients are generalized Vandermonde determinants).
This proves \reff{eq.exam.FH}.

\section*{Acknowledgments}

The research of SF was supported in part by an NSERC Discovery grant
(RGPIN-2014-06036).
The research of ADS was supported in part by EPSRC grant EP/N025636/1.


\clearpage

\addtocounter{section}{1}
\section*{Corrigendum}

Apoorva Khare has kindly drawn our attention to an error
in the proof of Theorem~3.2 in our paper \cite{FJS_totalpos}.
This theorem reads as follows:

\bigskip

\noindent
{\bf Theorem 3.2.}
\quad\hfill\vspace*{-1mm}
\begin{itemize}
   \item[(a)] {\em  Let $A$ be a Hankel matrix, all of whose contiguous minors
      of size $\le r$ are nonnegative.  Then $A$ is TN${}_r$.}
   \item[(b)] {\em  Let $A$ and $B$ be Hankel matrices,
      all of whose contiguous minors of size $\le r$ are nonnegative.
      Then $A+B$ is TN${}_r$.}
   \item[(c)] {\em  Let $A$ (resp.\ $B$) be a Hankel matrix,
      all of whose contiguous minors of size $\le r$ are nonnegative
      (resp.\ positive).
      Then $A+B$ is TP${}_r$.}
\end{itemize}

\medskip

\noindent
We attempted to prove (c) by invoking Sylvester's criterion
\cite[Theorem~7.2.5(a)]{Horn_13bis}
\linebreak
to conclude that every contiguous submatrix of $A$ of size $\le r$
is positive-semidefinite.
But the hypothesis of \cite[Theorem~7.2.5(a)]{Horn_13bis}
is that every {\em principal}\/ minor of the given
real symmetric (or complex hermitian) matrix is nonnegative;
it is {\em not}\/ sufficient to assume that every {\em contiguous}\/ minor
is nonnegative.
For instance,
$\displaystyle{ \begin{bmatrix}
                       0 & 0 & 1 \\
                       0 & 0 & 0 \\
                       1 & 0 & 0
                \end{bmatrix}
              }$
is a real symmetric matrix with all contiguous minors nonnegative,
but it is {\em not}\/ positive-semidefinite (in~fact, it has inertia +$-$0).
Similarly,
$\displaystyle{ \begin{bmatrix}
                       0 & 0 & 1 & 0 \\
                       0 & 0 & 0 & 1 \\
                       1 & 0 & 0 & 0 \\
                       0 & 1 & 0 & 0
                \end{bmatrix}
              }$
is a {\em nonsingular}\/ real symmetric matrix
with all contiguous minors nonnegative,
but it also fails to be positive-semidefinite (it has inertia ++$-${}$-$).
So our proof of Theorem~3.2(c) was fatally flawed.
Moreover, our proof of Theorem~3.2(b) was not as clear as it should have been.

Nevertheless, Theorem~3.2 is true ---
which is fortunate, since many later results in our paper depend on it
(notably Corollaries~3.3--3.5 and Theorems~4.5 and 5.11).
Here we would like to give a simple proof of Theorem~3.2
that follows the same strategy as our erroneous proof in \cite{FJS_totalpos},
but repairs the error.
The following easy lemma will be useful:

\begin{lemma}
  \label{lemma.hankelfix}
Let $A$ be a Hankel matrix.  Then the following are equivalent:
\begin{itemize}
   \item[(a)] Every minor of $A$ of size $\le r$ is positive
      (that is, $A$ is TP${}_r$).
   \item[(b)] Every contiguous minor of $A$ of size $\le r$ is positive.
   \item[(c)] Every contiguous submatrix of $A$ of size $\le r$ is
      positive-definite.
\end{itemize}
\end{lemma}

\proof
(a)$\implies$(b) is trivial,
and (b)$\implies$(a) is \cite[Lemma~2.4]{FJS_totalpos}.
(These results hold for arbitrary matrices $A$, not just Hankel matrices.)

(b)$\implies$(c):  Let $A'$ be a contiguous submatrix of $A$ of size $\le r$.
Since $A$ is Hankel, $A'$ is symmetric \cite[Lemma~2.7]{FJS_totalpos}.
And since every leading principal minor of $A'$ is positive,
$A'$ is positive-definite by Sylvester's criterion
\cite[Theorem~7.2.5(b)]{Horn_13bis}.

(c)$\implies$(b):  Let $A'$ be a contiguous submatrix of $A$ of size $\le r$.
If $A'$ is positive-definite, then $\det A' > 0$.
\qed

\proofof{Theorem~3.2}
(c) Let $A$ be a Hankel matrix, all of whose contiguous minors of size $\le r$
are nonnegative;
and let $B$ be a Hankel matrix, all of whose contiguous minors of size $\le r$ 
are positive.
We will prove that, for every $1 \le k \le r$,
every $k \times k$ minor of $A + \epsilon B$ is positive for all $\epsilon > 0$;
this implies in particular that $A + B$ is TP${}_r$.
The proof is by induction on $k$.
For $k=1$ the result is obvious (and does not need $A$ and $B$ to be Hankel).
So let $2 \le k \le r$, and assume that, for every $\epsilon > 0$,
every minor of $A + \epsilon B$ of size $\le k-1$ is positive.
Then, taking $\epsilon \to 0$,
we conclude that every minor of $A$ of size $\le k-1$ is nonnegative.
Now consider any $k \times k$ contiguous submatrix of $A$, call it $A'$;
and let $B'$ be the corresponding submatrix of $B$.
We have just shown that every minor of $A'$ of size $\le k-1$
is nonnegative;  and by hypothesis $\det A'$ is nonnegative.
So $A'$ is totally nonnegative.
In particular, since $A'$ is symmetric
(here is where we use the hypothesis that $A$ is Hankel
 \cite[Lemma~2.7]{FJS_totalpos})
and all its principal minors are nonnegative,
$A'$ is positive-semidefinite \cite[Theorem~7.2.5(a)]{Horn_13bis}.
Likewise, $B'$ is symmetric (since $B$ is Hankel)
and it is positive-definite by Lemma~\ref{lemma.hankelfix}(b)$\implies$(c).
Therefore $A' + \epsilon B'$ is positive-definite for all $\epsilon > 0$;
in particular, $\det(A' + \epsilon B') > 0$ for all $\epsilon > 0$.
This shows that all the $k$-by-$k$ contiguous minors of $A + \epsilon B$
are positive, for all $\epsilon > 0$.
Combining this with the inductive hypothesis
and using Lemma~\ref{lemma.hankelfix}(b)$\implies$(a),
we conclude that every minor of $A + \epsilon B$ of size $\le k$ is positive,
for all $\epsilon > 0$.
This completes the inductive step and hence proves~(c).

(a) The proof given in \cite{FJS_totalpos} is correct,
but for completeness we repeat it here.
Given $A$, we can take $B$ to be any TP Hankel matrix:
specifically, by \cite[Theorem~2.8]{FJS_totalpos}
we can take $B$ to be the Hankel matrix associated to
any Stieltjes moment sequence of infinite support
(for instance, $a_k = k!$ or $a_k = \lambda^{k^2}$ with $\lambda > 1$).
Then, by~(c), $A+\epsilon B$ is TP${}_r$ for all $\epsilon > 0$;
taking $\epsilon \to 0$ we conclude that $A$ is TN${}_r$.
This proves~(a).

(b) Finally, let $A$ and $B$ be Hankel matrices,
all of whose contiguous minors of size $\le r$ are nonnegative.
Then, by (a), $A$ and $B$ are TN${}_r$;
in particular, every contiguous submatrix of size $\le r$ of $A$ or $B$ is
positive-semidefinite by \cite[Theorem~7.2.5(a)]{Horn_13bis}
because all its principal minors are nonnegative.
So every contiguous submatrix of size $\le r$ of $A+B$
is positive-semidefinite,
and in particular has a nonnegative determinant.
Applying (a) to $A+B$, we obtain~(b).
\qed

We can now deduce an analogue of Lemma~\ref{lemma.hankelfix}
where positivity is replaced by nonnegativity:

\begin{corollary}
  \label{cor.hankelfix}
Let $A$ be a Hankel matrix.  Then the following are equivalent:
\begin{itemize}
   \item[(a)] Every minor of $A$ of size $\le r$ is nonnegative
      (that is, $A$ is TN${}_r$).
   \item[(b)] Every contiguous minor of $A$ of size $\le r$ is nonnegative.
   \item[(c)] Every contiguous submatrix of $A$ of size $\le r$ is
      positive-semidefinite.
\end{itemize}
\end{corollary}

\proof
(a)$\implies$(b) is trivial,
and (b)$\implies$(a) is Theorem~3.2(a).

(a)$\implies$(c):  Let $A'$ be a contiguous submatrix of $A$ of size $\le r$.
Since $A$ is Hankel, $A'$ is symmetric \cite[Lemma~2.7]{FJS_totalpos}.
And since every principal minor of $A'$ is nonnegative,
$A'$ is positive-semidefinite by Sylvester's criterion
\cite[Theorem~7.2.5(a)]{Horn_13bis}.

(c)$\implies$(b):  Let $A'$ be a contiguous submatrix of $A$ of size $\le r$.
If $A'$ is positive-semidefinite, then $\det A' \ge 0$.
\qed

\section*{Acknowledgments}

We are extremely grateful to Apoorva Khare for pointing out the mistake
in our proof of Theorem~3.2,
and to an anonymous referee for very helpful comments on an earlier draft
of this corrigendum.
We are also grateful to Alex Dyachenko for showing us a different correct proof
of Theorem~3.2, which is more complicated than the one given here
but which gives additional information about the vanishing and nonvanishing
of the minors of $A$.
We also thank Lily Li Liu for helpful discussions.

The research of SF was supported in part by an NSERC Discovery grant
(RGPIN-2014-06036).
The research of ADS was supported in part by EPSRC grant EP/N025636/1.


\addcontentsline{toc}{section}{References}
\begin{thebibliography}{99}


\bibitem{Akhiezer_65}  N.I. Akhiezer, {\em The Classical Moment Problem
   and Some Related Questions in Analysis}\/,
   translated by N.~Kemmer
   (Hafner, New York, 1965).

\bibitem{Ando_87}  T. Ando, Totally positive matrices,
   Lin. Alg. Appl. {\bf 90}, 165--219 (1987).

%

\bibitem{Berman_03}  A. Berman and N. Shaked-Monderer,
   {\em Completely Positive Matrices}\/
   (World Scientific, River Edge NJ, 2003).

%



\bibitem{Crans_01}  A.S. Crans, S.M. Fallat and C.R. Johnson,
   The Hadamard core of the totally nonnegative matrices,
   Lin. Alg. Appl. {\bf 328}, 203--222 (2001).

\bibitem{Cryer_73}  C.W. Cryer, The {\em LU}\/-factorization of totally
   positive matrices, Lin. Alg. Appl. {\bf 7}, 83--92 (1973).

\bibitem{Cryer_76}  C.W. Cryer, Some properties of totally positive matrices,
   Lin. Alg. Appl. {\bf 15}, 1--25 (1976).

\bibitem{Curto_91}  R. Curto and L.A. Fialkow,
   Recursiveness, positivity, and truncated moment problems,
   Houston J. Math. {\bf 17}, 603--635 (1991).

\bibitem{Fallat_07}  S.M. Fallat and C.R. Johnson,
   Hadamard powers and totally positive matrices,
   Lin. Alg. Appl. {\bf 423}, 420--427 (2007).

\bibitem{Fallat_11}  S.M. Fallat and C.R. Johnson,
   {\em Totally Nonnegative Matrices}\/
   (Princeton University Press, Princeton NJ, 2011).

\bibitem{Fekete_12}  M. Fekete and G. P\'olya,  \"Uber ein Problem von
   Laguerre, Rendiconti Circ. Mat. Palermo {\bf 34}, 89--120 (1912).

\bibitem{Fitzgerald_77}  C.H. FitzGerald and R.A. Horn,
   On fractional Hadamard powers of positive definite matrices,
   J. Math. Anal. Appl. {\bf 61}, 633--642 (1977).

\bibitem{Fomin_00}  S. Fomin and A. Zelevinsky,
   Total positivity: tests and parametrizations,
   Math. Intelligencer {\bf 22}, no.~1, 23--33 (2000).

\bibitem{Gantmacher_59}  F.R. Gantmacher, {\em The Theory of Matrices}\/,
   2 vols.\ (Chelsea, New York, 1959).
   Reprinted by AMS Chelsea, Providence, RI, 1998.

\bibitem{Gantmacher_02}  F.R. Gantmacher and M.G. Krein,
   {\em Oscillation Matrices and Kernels and Small Vibrations of
       Mechanical Systems}\/
   (AMS Chelsea Publishing, Providence RI, 2002).
   Based on the second Russian edition, 1950.

\bibitem{Gantmakher_37}  F. Gantmakher and M. Krein,
   Sur les matrices compl\`etement non n\'egatives et oscillatoires,
   Compositio Math. {\bf 4}, 445--476 (1937).

\bibitem{Gasca_96}  M. Gasca and C.A. Micchelli, eds.,
   {\em Total Positivity and its Applications}\/
   (Kluwer, Dordrecht, 1996).

\bibitem{Gasca_92}  M. Gasca and J.M. Pe\~na,
   Total positivity and Neville elimination,
   Lin. Alg. Appl. {\bf 165}, 25--44 (1992).

\bibitem{Holtz_12}  O. Holtz and M. Tyaglov,
   Structured matrices, continued fractions, and root localization of
   polynomials, SIAM Rev. {\bf 54}, 421--509 (2012).

\bibitem{Horn_91}  R.A. Horn and C.R. Johnson,
   {\em Topics in Matrix Analysis}\/
   (Cambridge University Press, Cambridge, 1991).

\bibitem{Horn_13}  R.A. Horn and C.R. Johnson,
   {\em Matrix Analysis}\/, 2nd edition
   (Cambridge University Press, Cambridge, 2013).

\bibitem{Iohvidov_82}  I.S. Iohvidov, {\em Hankel and Toeplitz Matrices and
   Forms}\/ (Birkh\"auser, Boston, Mass., 1982).

\bibitem{Jameson_06}  G.J.O. Jameson, Counting zeros of generalized
   polynomials: Descartes' rule of signs and Laguerre's extensions,
   Math. Gazette {\bf 90}, 223--234 (2006).

\bibitem{Johnson_12}  C.R. Johnson and O. Walch,
   Critical exponents: Old and new,
   Electronic J. Lin. Alg. {\bf 25}, 72--83 (2012).

\bibitem{Karlin_68}  S. Karlin, {\em Total Positivity}\/
   (Stanford University Press, Stanford CA, 1968).

\bibitem{Laguerre_1883}  E. Laguerre, Sur la th\'eorie des \'equations
   num\'eriques, Journal de Math\'ematiques Pures et Appliqu\'ees
   (3) {\bf 9}, 99--146 (1883).
   Also in {\em \OE{}uvres de Laguerre}\/, tome I, pp.~3--47.

\bibitem{Markham_70}  T.L. Markham, A semigroup of totally nonnegative
   matrices, Lin. Alg. Appl. {\bf 3}, 157--164 (1970).


\bibitem{Pinkus_10}  A. Pinkus, {\em Totally Positive Matrices}\/
   (Cambridge University Press, Cambridge, 2010).

\bibitem{Polya-Szego}  G. P\'olya and G. Szeg\H{o},
   {\em Problems and Theorems in Analysis II}\/
   (Springer-Verlag, Berlin--Heidelberg--New York, 1976).

\bibitem{Rahman_02}  Q.I. Rahman and G. Schmeisser,
   {\em Analytic Theory of Polynomials}\/
   (Clarendon Press, Oxford, 2002).

\bibitem{Schoenberg_30}  I. Schoenberg, 
   \"Uber variationsvermindernde lineare Transformationen,
   Math. Z. {\bf 32}, 321--328 (1930).

\bibitem{Shohat_43}  J.A. Shohat and J.D. Tamarkin,
   {\em The Problem of Moments}\/
   (American Mathematical Society, New York, 1943).
   Revised edition, 1950.

\bibitem{Sokal_totalpos}  A.D. Sokal,
   Coefficientwise total positivity (via continued fractions)
   for some Hankel matrices of combinatorial polynomials,
   in preparation.




\end{thebibliography}

\addcontentsline{toc}{section}{References}
\begin{thebibliography}{99}

\bibitem{FJS_totalpos}  S.M. Fallat, C.R. Johnson and A.D. Sokal,
   Total positivity of sums, Hadamard products and Hadamard powers:
   Results and counterexamples,
   Lin. Alg. Appl. {\bf 520}, 242--259 (2017).

\bibitem{Horn_13bis}  R.A. Horn and C.R. Johnson,
   {\em Matrix Analysis}\/, 2nd edition
   (Cambridge University Press, Cambridge, 2013).


\end{thebibliography}
\end{document}